\documentclass[11pt,twoside]{article}

\date{}
\usepackage{amsmath}
\usepackage{amssymb}
\usepackage{amsthm}

\newcommand{\C}{\mathbb{C}}
\newcommand{\R}{\mathbb{R}}
\newcommand{\g}{\mathfrak{g}}

\newcommand{\lk}{\mathfrak{k}}

\newcommand{\bdm}{\begin{displaymath}}
\newcommand{\edm}{\end{displaymath}}

\theoremstyle{definition}

\newtheorem{lem}{Lemma}
\newtheorem{thm}{Theorem}
\newtheorem{defn}{Definition}
\newtheorem{prop}{Proposition}
\newtheorem{rem}{Remark}

\newtheorem{axiom}{Corollary}

\title{{Bianchi-B\"{a}cklund transforms
and dressing actions, revisited}}
\author{R. Pacheco}
\begin{document}
\maketitle \emph{Departamento de Matem\'{a}tica, Universidade da
Beira Interior, Rua Marqu\^{e}s d'\'{A}vila e Bolama 6201-001
Covilh\~{a} - Portugal.} \emph{Email: rpacheco@mat.ubi.pt}

%
%

\section{Introduction} The Bianchi-B\"{a}cklund transforms, which were introduced
by Bianchi \cite{Bi} for surfaces of positive constant Gauss
curvature (equivalently, solutions of $\sinh$-Gordon equation),
are an extension of the B\"{a}cklund transforms \cite{Ei} for
surfaces of negative constant Gauss curvature (equivalently,
solutions of sine-Gordon equation). In contrast to the negative
case, a Bianchi-B\"{a}cklund transform of a real surface is in
general complex. In order to obtain a new real surface with
positive constant Gauss curvature (CGC $K>0$), one has to apply
two successive Bianchi-B\"{a}cklund transforms, where the second
transform has to be matched to the first in a particular way.

In \cite{Ma}, A. Mahler showed in a constructive way that the
classical Bianchi-B\"{a}cklund procedure for obtaining  a new real
 CGC $K>0$ surface
$\tilde{f}$ out of an old one $f$ amounts to dressing the extended
framing associated to the Gauss map of $f$, which is an harmonic
map, by a certain dressing matrix. In this paper we shall give an
alternative approach to Mahler's work. Following the philosophy
developed by Terng and Uhlenbeck \cite{Tu}, we start with certain
basic elements, the {\textit{simple factors}}, for which the
dressing action can be computed explicitly. We show that each
single Bianchi-B\"{a}cklund transform  corresponds to the dressing
action of a certain simple factor. A nice geometrical
parameterization of these simple factors is available and we shall
see how to relate it with the classical parameterization of
Bianchi-B\"{a}cklund transforms. As a consequence, we recover the
result announced by Mahler. Moreover:

Sophus Lie proved that every  B\"{a}cklund transformation is a
combination of transformations of Lie and Bianchi. Lie
transformations correspond to the invariance of sine-Gordon
equation under Lorentz transforms and Bianchi transformations are
 B\"{a}cklund transformations for which the tangent spaces at corresponding points on the original surface and
 the new surface are orthogonal. On the other hand, Bonnet observed that
 $\sinh$-Gordon equation admits an invariance similar to that of Lie.
In this paper, we shall show how this invariance of
$\sinh$-Gordon equation can be combined with those
Bianchi-B\"{a}cklund transformations satisfying the orthogonality
condition of Bianchi in order to produce all Bianchi-B\"{a}cklund
transformations.

Finally, we shall mention that this correspondence between
Bianchi-B\"{a}cklund transforms and dressing actions of simple
type has already been established within another setting
\cite{Buiso,HP,IK}. In fact, each CGC $K>0$ surface corresponds to
a pair of parallel constant mean curvature (CMC) surfaces. Hence,
by applying the two-step Bianchi-B\"{a}cklund procedure, we can
transform a CMC surface in a new CMC surface. Since CMC surfaces
are isothermic, they also allow Darboux transformations via sphere
congruences. Burstall \cite{Buiso} showed that Darboux
transformations for isothermic surfaces correspond to certain loop
group actions of the simple type. On the other hand, Darboux
transforms of CMC surfaces are equivalent to Bianchi-B\"{a}cklund
transforms of CMC surfaces \cite{HP,IK}. Hence, the two-step
Bianchi-B\"{a}cklund procedure for CGC $K>0$ surfaces corresponds
to
 loop group actions of the simple type. In our  approach, we study CGC $K>0$ surfaces via their Gauss
 maps, which are harmonic maps into the unit sphere $S^2$. Since
 there is no well established theory of CMC surfaces unifying
 the harmonic map and isothermic surface theories, we do not know
 how to relate our dressing action to that of F. Burstall
 \cite{Buiso} (the underlying symmetry groups seem
 quite different).

\emph{Acknowledgments.} I would like to express my deepest
gratitude to F. Burstall for calling my attention to this topic
while I was a PhD candidate at University of Bath.

\section{CMC Surfaces, CGC Surfaces, and Harmonic Maps}
We start by reviewing some well known facts about CMC and CGC
surfaces, and their correspondence with harmonic maps into $S^2$.
We refer the reader to \cite{FH} for details.

\vspace{.20in}

Consider on $\mathbb{C}$ the standard conformal structure and let
$(\cdot,\cdot)$ be the standard inner product of $\mathbb{R}^3$. A
smooth map $\varphi:\C\rightarrow S^2$ is harmonic if, and only
if, the one form $\varphi\times *{{d}}\varphi$ is closed. Hence,
if $\varphi:\C\rightarrow S^2$ is harmonic, there exists
$F:\C\rightarrow \mathbb{R}^3$ such that ${{d}}F=\varphi\times
*{{d}}\varphi$. Moreover, $F$ is an immersion if, and only if,
$\varphi$ is an immersion. In this case: $F$ is a CGC $K=1$
surface; the conformal structure on $\C$ is given by the second
fundamental form $\mathrm{\Pi}_F$ of $F$; and  $F$ has an umbilic
point at $p\in\C$ if, and only if, $\varphi$ is conformal at $p$.
Conversely, given  a CGC $K=1$ immersion $F:\C \rightarrow
\mathbb{R}^3$, $\mathrm{\Pi}_F$ is definite and the corresponding
Gauss map $\varphi:(\C,\mathrm{\Pi}_F) \rightarrow S^2$ is
harmonic.

Any harmonic map $\varphi:\C\to S^2$ which is everywhere
non-conformal is the Gauss map of two conformal CMC $H=\frac12$
immersions. These immersions have no
 umbilic points and are given by $f^\pm=F\pm\varphi:\C \to
\mathbb{R}^3$, where $F$ (not necessarily an immersion) is such
that $dF=\varphi\times *{{d}}\varphi$. Conversely, if $f:\C \to
\mathbb{R}^3$ is a conformal CMC $H=\frac12$ immersion without
umbilic points, the corresponding Gauss map $\varphi:\C\to S^2$ is
such that $dF=\varphi\times *{{d}}\varphi$, with $F=f+\varphi$,
which means that $\varphi$ is harmonic and everywhere
non-conformal.

When a conformal CMC $H=\frac{1}{2}$ immersion has no umbilic
points, it is well known that we can always choose a conformal
coordinate $z=x+iy$ with respect to which the second fundamental
form of $f$ is diagonal, that is, $z$ is a conformal curvature
line coordinate on $\C$. More precisely, we have
\begin{equation*}
\mathrm{I}_f=e^{2\omega}({d}x^2+{d}y^2),\,\,
\mathrm{\Pi}_f=e^\omega(\sinh\omega\,{d}x^2+ \cosh \omega\,
{d}y^2)\,,
\end{equation*}
where $\omega:\C\rightarrow\mathbb{R}$ is a solution to the Gauss
($\sinh$-Gordon) equation
\begin{equation}\label{sinhgordon}
\triangle\omega+\sinh \omega \cosh \omega=0\,.
\end{equation}
So, away from the points where $\omega$ vanishes, $F=f+\varphi$ is
a CGC $K=1$ immersion and we have
\begin{equation}\label{firstandsecond}
\mathrm{I}_F=\cosh^2\omega
\,{d}x^2+\sinh^2\omega\,{d}y^2,\,\,\mathrm{\Pi}_F=-\sinh\omega\cosh\omega
({d}x^2+{d}y^2)\,.
\end{equation}

Bonnet observed that if $\omega$ is a solution of
(\ref{sinhgordon}) and $\sigma\in\C$, then $\omega\circ R_\sigma$
is also a solution, where
$$R_\sigma(x,y)=(\cos \sigma \,x-\sin\sigma \,y, \sin \sigma \,x +\cos\sigma \,y)$$ (cf. \cite{Bi}). Hence, to every CGC
$K=1$ surface $F$, there is an associated $S^1$-family of CGC
$K=1$ surfaces. We can enlarge this family: if $\omega$ is a
solution of $\sinh$-Gordon equation, then so is
$\omega^\lambda(z,\bar{z})=\omega(\lambda z,\lambda^{-1}\bar{z})$,
for each $\lambda\in\mathbb{C}^*$; hence, there is an associated
smooth $\C^*$-family of CGC $K=1$ surfaces $F^\lambda$. We denote
$F^\lambda=\mathbf{S}_\sigma(F)$, with $\lambda=e^{-i\sigma}$. In
general, for $\lambda \notin S^1$, observe that $\omega^\lambda$
is a complex solution and, consequently, $F^\lambda$ is a complex
surface.

\vspace{.20in}

 \textit{Pseudospherical surfaces.}
Denote by $(\mathbb{R}^2,\epsilon)$ the Minkowski space, where
$\epsilon:=dudv$ is the Minkowski metric in characteristic
coordinates: $u=x+y$ and $v=x-y$. Define the star operator on
one-forms by $*(Adu+Bdv)=Adu-Bdv$. Again,  a smooth map
$\varphi:(\mathbb{R}^2,\epsilon)\rightarrow S^2$ is (Lorentz)
harmonic if, and only if, the one form $\varphi\times
*{{d}}\varphi$ is closed (cf. \cite{MS}). Hence, if
$\varphi:(\mathbb{R}^2,\epsilon)\rightarrow S^2$ is harmonic,
 there exists
$F:\mathbb{R}^2\rightarrow \mathbb{R}^3$ such that
${{d}}F=\varphi\times *{{d}}\varphi$. Again, $F$ is an immersion
if, and only if, $\varphi$ is an immersion. In this case: $F$ is a
CGC $K=-1$ surface; the Minkowski metric $\epsilon$ and the second
fundamental form  $\mathrm{\Pi}_F$ are conformally equivalent; if
$\varphi$ is \emph{weakly regular}, that is, if $d\varphi$ never
vanishes on the characteristic directions, the corresponding
surface admits a global parametrization in asymptotic coordinates
$\xi,\eta$ such that
 the two fundamental forms become
\begin{equation*}
\mathrm{I}_F={d}\xi^2+2\cos \omega \, d\xi d\eta+{d}\eta^2,\,\,
\mathrm{\Pi}_F=2\sin \omega\,{d}\xi d\eta\,,
\end{equation*}
where $\omega$ satisfies the \emph{sine-Gordon} equation
$\omega_{\xi\eta}=\sin\omega$ (cf. \cite{MS}).  Lie observed that
this equation is invariant under Lorentz transformations:  if
$\omega$ is a solution, then so is
$\omega^\lambda(\xi,\eta)=\omega (\lambda\xi, \lambda^{-1} \eta)$,
with $\lambda$ belonging to the multiplicative group
$\mathbb{R}^*$ of nonzero real numbers. Hence, to every CGC $K=-1$
surface $F$, we can associate a smooth $\mathbb{R}^*$-family of
CGC $K=-1$ surfaces $F^\lambda$. $F^\lambda$ is called a \emph{Lie
transform} of $F$ and we denote $F^\lambda=\mathbf{L}_\lambda(F)$.

\section{Bianchi-B\"{a}cklund Transforms}

The Bianchi-B\"{a}cklund transforms, which were introduced by
Bianchi \cite{Bi} for positive CGC surfaces, are an extension of
the B\"{a}cklund transformation \cite{Ei} for surfaces of negative
constant Gauss curvature. We shall now review briefly this theory,
following \cite{SW}.

\vspace{.25in}

Let $F:\mathbb{C}\rightarrow\mathbb{R}^3$ be a CGC $K=1$ surface
 and $z=x+iy$ a curvature line coordinate on
$\C$, with fundamental forms $\mathrm{I}_F$ and
$\mathrm{\Pi}_F$ given by (\ref{firstandsecond}), where $\omega$
is a solution to the $\sinh$-Gordon equation (\ref{sinhgordon}).
In particular, $F$ has no umbilic points and the coordinate $z$ is
conformal for the parallel CMC surface but not for $F$. Let
$\varphi$ be the corresponding Gauss map. Define the orthonormal
frame $[e_1,e_2,e_3]$ by:
\begin{equation*}
e_1  =  \frac{1}{\cosh\omega}F_x\,,\,\,\,  e_2  =
\frac{1}{\sinh\omega}F_y\,,\,\,\,e_3  = e_1\times e_2 \,.
\end{equation*}
Let $\tilde{F}:\C\rightarrow (\mathbb{R}^3)^\C\cong \C^3$ be a
complex  surface   with Gauss map $\tilde{\varphi}$.
\begin{defn}\cite{SW}
\emph{We say that $\tilde{F}$ is a {Bianchi-B\"{a}cklund
transform} of the CGC $K=1$ surface $F$ if it satisfies the
following properties:}  \emph{$z$ is a curvature line coordinate
with respect to $F$ and $\tilde{F}$;}
\emph{$(\tilde{F}-F,\varphi)=(\tilde{F}-F,\tilde{\varphi})=0$;
$\tilde{F}-F$ has constant length; }  \emph{the normals have a
constant angle  with each other. }
\end{defn}

So, suppose that $\tilde{F}$ is a Bianchi-B\"{a}cklund transform
of $F$. Let $\phi:\mathbb{C}\rightarrow \C$ be the angle formed by
the tangent line $\tilde{F}-F$ and $e_1$. Then,
$\tilde{F}=F+\mu(\cos \phi e_1+ \sin\phi e_2)$ for some
$\mu\in\C\setminus\{0\}$. We also have
\begin{equation}\label{osigma}
(\varphi,\tilde{\varphi})=\cos\sigma\,\,\,\,\,\,\,\,\,\,\,\,\mbox{and}\,\,\,\,\,\,\,\,\,\,\,\,\varphi\times
\tilde{\varphi}=\sin \sigma( \cos \phi e_1+ \sin \phi e_2 ),
\end{equation}
for some constant angle $\sigma$.
\begin{thm}\cite{SW}
\emph{$\tilde{F}$ is a new CGC $K=1$ surface and
\begin{equation}\label{newcgc}
\tilde{F}=F+\frac{1}{\sinh\beta}(\cosh \theta e_1+i \sinh\theta
e_2)\,,
\end{equation}
where $\beta\in\C\setminus\{in\pi, n\in\mathbb{Z}\}$ is a
constant, $\mu=\frac{1}{\sinh\beta}$, $\cot \sigma=-i\cosh \beta$,
and $\theta:\mathbb{C}\rightarrow\C$, with $\theta =-i\phi$, is a
solution of
\begin{equation}\label{bbpde}
\left\{\begin{array}{ll} \theta_x+i\omega_y &= \sinh \beta \sinh
\theta \cosh \omega+\cosh \beta \cosh \theta \sinh \omega
\\ i\theta_y+\omega_x &=-\sinh \beta\cosh\theta\sinh\omega-
\cosh \beta \sinh \theta\cosh \omega \end{array}\right.\,.
\end{equation}
The first and second fundamental forms of $\tilde{F}$ are given by
\begin{equation*}
\mathrm{I}_{\tilde{F}}=\cosh^2\theta
{d}x^2+\sinh^2\theta{d}y^2,\,\,\mathrm{\Pi}_{\tilde{F}}=-\sinh\theta\cosh\theta
({d}x^2+{d}y^2)\,.
\end{equation*}}
\end{thm}
Hence, the classical  Bianchi-B\"{a}cklund transformations are
determined by two complex numbers: the spectral parameter
$\beta\in \C\setminus\{in\pi, n\in\mathbb{Z}\}$ and an initial
angle $\theta_0\in\C$.

From the analytical point of view, the Bianchi-B\"{a}cklund
transformations may be interpreted  as a procedure of obtaining
new solutions of the $\sinh$-Gordon equation from an old one. In
fact, if $\omega:\mathbb{R}^2\rightarrow \C$ is a solution of the
$\sinh$-Gordon equation, then any solution $\theta$ of the
Bianchi-B\"{a}cklund PDEs (\ref{bbpde}) will also solve the
$\sinh$-Gordon equation.

We denote by ${{\mathbf{BB}}}_{\beta}(F)$  the family of
Bianchi-B\"{a}cklund transforms of $F$  with spectral parameter
$\beta$ and by ${{\mathbf{BB}}}_{\beta,\theta_0}(F)$
 the
Bianchi-B\"{a}cklund transform of $F$  determined by
$(\beta,\theta_0)$.
\begin{lem}\label{-beta}
\emph{Suppose that $\beta_1,\beta_2\in\C\setminus\{in\pi,
n\in\mathbb{Z}\}$ satisfy $\sinh\beta_1=-\sinh\beta_2$ and
$\cosh\beta_1=-\cosh\beta_2$. Then
${\mathbf{BB}}_{\beta_1}(F)={\mathbf{BB}}_{\beta_2}(F)$.}
\end{lem}
\begin{proof}
 Let $\theta_1$ be a solution of (\ref{bbpde}) for $\beta_1$.
Then $\theta_2=i\pi+\theta_1$ is a solution of (\ref{bbpde}) for
$\beta_2$. Taking account formula (\ref{newcgc}), we conclude that
the pairs $(\beta_1,\theta_1)$ and $(\beta_2,\theta_2)$ produce
the same surface, and we are done.
\end{proof}

In contrast to the negative CGC case, the solution $\theta$ of the
Bianchi-B\"{a}cklund PDEs is in general complex and so $\tilde{F}$
will be complex. To obtain a new real solution of $\sinh$-Gordon
equation we must perform two iterations of  Bianchi-B\"{a}cklund
transformations:
\begin{thm}\label{realcondb}\cite{SW}
\emph{Start with the real CGC $K=1$ surface $F$. If the reality
condition $\beta^*=i\pi-\overline{\beta}$  holds (hence
$\theta^*=-\overline{\theta}$),  then
${\mathbf{BB}}_{\beta^*,\theta^*_0}\big({\mathbf{BB}}_{\beta,\theta_0}(F)\big)$
is real.}
\end{thm}
Finally  we recall the Bianchi--B\"{a}cklund Permutability
theorem:
\begin{thm}\cite{SW}
\emph{Let $F$ be a CGC $K=1$ surface and
$\beta,\beta^*\in\C\setminus\{in\pi, n\in\mathbb{Z}\}$. Then
$${\mathbf{BB}}_{\beta^*}\big({\mathbf{BB}}_{\beta}(F)\big)=
{\mathbf{BB}}_{\beta}\big({\mathbf{BB}}_{\beta^*}(F)\big).$$}
\end{thm}

\vspace{.20in}

 \textit{B\"{a}cklund transforms.} Let $F:\mathbb{R}^2\to \mathbb{R}^3$ be a $K=-1$ pseudospherical  surface.
A  surface $\tilde{F}:\mathbb{R}^2\to \mathbb{R}^3$ is a
\emph{B\"{a}cklund transform} of $F$ if it satisfies the following
properties \cite{Ei}:  The coordinates $u,v$ correspond to
parametrization along asymptotic lines with respect to $F$ and
$\tilde{F}$;
$(\tilde{F}-F,\varphi)=(\tilde{F}-F,\tilde{\varphi})=0$;
$\tilde{F}-F$ has constant length; the normals have a constant
angle $\sigma$ with each other. We denote
$\tilde{F}=\mathbf{B}_\beta(F)$, with $\beta=\tan \sigma/2$. In
the particular case $\beta=1$ (the tangent planes at corresponding
points on $F$ and $\tilde{F}$ are orthogonal), $\tilde{F}$ is
called a \emph{Bianchi transform} of $F$. Lie observed that every
B\"{a}cklund transformation is a combination of transformations of
Lie and Bianchi (cf. \cite{Ei}). More precisely,
$\mathbf{B}_\beta=\mathbf{L}_\beta^{-1}\circ\mathbf{B}_{1}\circ\mathbf{L}_\beta$.

\section{Harmonic Maps and Loop Groups}
In this section we review the reformulation of harmonicity
equations, for maps from a Riemann surface into a compact
symmetric space, in terms of loops of flat connections, referring
the reader to \cite{BP,U} for details.

\vspace{.20in}

Let $G$ be a compact (connected) semisimple  matrix Lie group,
with identity $e$ and Lie algebra $\mathfrak{g}$. Equip $G$ with a
bi-invariant metric. Let $G^{\mathbb{C}}$ be the complexification
of $G$, with Lie algebra $\mathfrak{g}^{\mathbb{C}}$ (thus
$\mathfrak{g}^\mathbb{C}=\mathfrak{g}\otimes \mathbb{C}$).
Consider a symmetric space $N=G/K$ with automorphism $\tau$ and
associated symmetric decomposition $\g=\lk\oplus\mathfrak{m}$. Let
$\phi:\C\rightarrow N$ be a smooth map and take a lift
$\psi:\C\rightarrow G$ of $\phi$, that is, we have
$\phi=\pi\circ\psi$ where $\pi:G\rightarrow G/K$ is the coset
projection. Corresponding to the symmetric decomposition
$\g=\lk\oplus\mathfrak{m}$ there is a decomposition of
$\alpha=\psi^{-1}{d}\psi$, $
\alpha=\alpha_\mathfrak{k}+\alpha_\mathfrak{m}$.  Let
$\alpha_\mathfrak{m}=\alpha'_\mathfrak{m}+\alpha''_\mathfrak{m}$
be the type decomposition of $\alpha_\mathfrak{m}$ into
$(1,0)$-form and $(0,1)$-form of $\C$. Consider the loop of
$1$-forms
$\alpha_\lambda=\lambda^{-1}\alpha'_\mathfrak{m}+\alpha_\mathfrak{k}
+\lambda\alpha''_\mathfrak{m}$. We may view $\alpha_\lambda$ as a
$\Lambda_\tau\g$-valued $1$-form, where
\begin{equation}
\label{lambdatau}
 \Lambda_\tau\g=\left\{\xi:S^1\rightarrow
\g\,\,(\mathrm{smooth})\,|\,\,
\tau(\xi(\lambda))=\xi(-\lambda)\,\,\mathrm{for\,\,
all}\,\,\lambda\in S^1\right\}.
\end{equation}
It is well known that $\phi$ is harmonic if, and only if,
${d}+\alpha_\lambda$ is a loop of flat connections on the trivial
bundle $\underline{\C}^n=\C\times\C^n$. Hence, if $\phi$ is
harmonic, we can define a smooth map $\Phi:\C\rightarrow
\Lambda_\tau G$, where $\Lambda_\tau G$ is the
infinite-dimensional Lie group corresponding to the loop Lie
algebra (\ref{lambdatau}),
\begin{equation*}
\Lambda_\tau G=\left\{\gamma:S^1\rightarrow
G\,\,(\mathrm{smooth})\,|\,\,
\tau(\gamma(\lambda))=\gamma(-\lambda)\,\,\mathrm{for\,\,all}\,\,\lambda\in
S^1\right\},
\end{equation*}
such that $\Phi^{-1}d\Phi=\alpha_\lambda$. The smooth map $\Phi$
is called an \textit{extended framing} (associated to $\phi$) and
gives rise to a smooth $S^1$-family of harmonic maps
$\phi^\lambda=\pi\circ \Phi_\lambda$ (here we are using the
notation $\Phi_\lambda(z)=\Phi (z)(\lambda)$), with $\phi^1=\phi$.

 Inspired by the classical theory of Bianchi-B\"{a}cklund transforms, where
we have to deal with complex surfaces,  we generalize the notion
of extended framing as follows: a smooth map $\Phi:\C\rightarrow
\Lambda_\tau G^\C$ is called a \textit{complex extended framing}
if $\Phi$ satisfies
$$\Phi^{-1}{d}\Phi=(\lambda^{-1}A+B){d}z+(C+\lambda D){d}\bar{z},$$ where
$B,C:\C\rightarrow \mathfrak{k}^\C$ and $A,D:\C\rightarrow
\mathfrak{m}^\C$. Of course, an extended framing is a complex
extended framing satisfying the reality condition
$\Phi(\lambda)=\overline{\Phi_{1/\bar{\lambda}}}$.

\vspace{.15in}

In order to recover our CGC $K=1$ surface and the corresponding
Gauss map from an extended framing, we proceed as follows:

 As $\mathrm{SO}(3,\mathbb{R})$-modules, we can identify
$(\mathbb{R}^3,\times)\cong\mathfrak{so}(3,\mathbb{R})$ via
$\mathbb{R}^3 \ni u\mapsto\xi_u\in\mathfrak{so}(3,\mathbb{R})$,
where $\xi_{u}(v)=u\times v$. The inner product on
$\mathfrak{so}(3,\mathbb{R})$ inherited from $\mathbb{R}^3$ is
$(\xi,\eta)=-\frac{1}{2}\mathrm{Tr}\,\xi\eta$. Let $e_1,e_2,e_3$
be the canonical orthonormal basis of $\mathbb{R}^3$ and
$K\subset\mathrm{SO}(3,\mathbb{R})$ the stabilizer of $e_1\in
S^2$. Consider the automorphism $\tau$ of
$\mathrm{SO}(3,\mathbb{R})$ given by conjugation by
\begin{equation}\label{Q}
Q=\begin{pmatrix}
  1 & 0 & 0 \\
  0 & {-1} & 0 \\
  0 & 0 & {-1}
\end{pmatrix}\,.
\end{equation}
Then $K$ is the identity component of  the fixed set $K^\tau$ of
$\tau$. The corresponding symmetric decomposition
$\mathfrak{so}(3,\mathbb{C})=\mathfrak{k}^\C\oplus
\mathfrak{m}^\C$ is given by:
\begin{equation}
\label{syme} \mathfrak{k}^\C=\{x\xi_{e_1}:\,x\in\C\},
\quad\quad\mathfrak{m}^\C=\{y\xi_{e_2}+z\xi_{e_3}:\,y,z\in\C\}.
\end{equation}

Given a complex extended framing $\Phi:\C\rightarrow \Lambda_\tau
\mathrm{SO}(3,\mathbb{C})$, consider the smooth map
$\varphi^\lambda:\C\rightarrow \mathfrak{so}(3,\mathbb{C})$ defined by
$\varphi^\lambda=\Phi_\lambda\cdot \varphi_0$, where $\varphi_0=\xi_{e_1}$. Then
one can easily prove that the smooth map
\begin{equation}\label{adlersymes}
F^\lambda=-i\lambda\frac{\partial \Phi}{\partial
\lambda}\Phi_\lambda^{-1}:\C\rightarrow
\mathfrak{so}(3,\mathbb{C})
\end{equation}
satisfies
${d}F^\lambda=\left[\varphi^\lambda,*{d}\varphi^\lambda\right].$
In particular, $F^\lambda$ is a (real or complex) CGC $K=1$
immersion if $\varphi^\lambda$ is an immersion. Moreover,
$F^\lambda=\mathbf{S}_\sigma(F),$ with $\lambda=e^{-i\sigma}$.

 \vspace{.20in}

\textit{Reformulation of Lorentz harmonicity in terms of loops of
flat connections.} Let $\varphi:(\mathbb{R}^2,\epsilon)\to S^2$ be
a smooth map and take a lift $\psi:\mathbb{R}^2\rightarrow
\mathrm{SO}(3,\mathbb{R})$ of $\varphi$, that is, we have
$\phi=\pi\circ\psi$ where
$\pi:\mathrm{SO}(3,\mathbb{R})\rightarrow S^2\simeq
\mathrm{SO}(3,\mathbb{R})/\mathrm{SO}(2,\mathbb{R})$ is the coset
projection. Corresponding to the  symmetric decomposition
(\ref{syme}) of $\mathfrak{so}(3,\mathbb{R})$,
$\mathfrak{so}(3,\mathbb{R})=\lk\oplus\mathfrak{m}$, there is a
decomposition of $\alpha=\psi^{-1}{d}\psi$, $
\alpha=\alpha_\mathfrak{k}+\alpha_\mathfrak{m}$. Write
$\alpha_\mathfrak{m}=A_1'du+A_1''dv$ and
$\alpha_\mathfrak{k}=A'_0du+A''_0dv$. Consider the loop of
$1$-forms
$$\alpha_\lambda=\lambda^{-1}A_1'du+A'_0du+A''_0dv
+\lambda A_1''dv,$$ where $u,v$ are characteristic coordinates,
with $\lambda$ belonging to the multiplicative group
$\mathbb{R}^*$ of nonzero real numbers. We may view
$\alpha_\lambda$ as a
$\Lambda^r_\tau\mathfrak{so}(3,\mathbb{R})$-valued $1$-form, where
$$\Lambda^r_\tau\mathfrak{so}(3,\mathbb{R})=\left\{\xi:\mathbb{C}^*\rightarrow
\mathfrak{so}(3,\mathbb{C})\,\,(\mathrm{smooth})\,|\,\,
\tau(\xi(\lambda))=\xi(-\lambda),\,\,
\overline{\xi(\lambda)}=\xi(\bar{\lambda}) \right\}.$$ Observe
that the reality condition
$\overline{\xi(\lambda)}=\xi(\bar{\lambda})$ forces any element
$\xi$  of $\Lambda^r_\tau\mathfrak{so}(3,\mathbb{R})$ to  assume
values in $\mathfrak{so}(3,\mathbb{R})$ for all
$\lambda\in\mathbb{R}^*$. The smooth map $\varphi$ is (Lorentz)
harmonic if, and only if, $d+\alpha_\lambda$ is a flat connection
for all $\lambda\in\mathbb{R}^*$  \cite{MS}. Hence, if $\varphi$
is harmonic, we can define a smooth map
$\Phi:\mathbb{R}^2\rightarrow \Lambda_\tau^r
\mathrm{SO}(3,\mathbb{R})$, where $\Lambda^r_\tau
\mathrm{SO}(3,\mathbb{R})$ is the infinite-dimensional Lie group
corresponding to $\Lambda^r_\tau\mathfrak{so}(3,\mathbb{R})$, such
that $\Phi^{-1}d\Phi=\alpha_\lambda$. Again, the smooth map $\Phi$
is called an \textit{extended framing}, $\varphi$ is recovered
from $\Phi$ via $\varphi=\pi\circ \Phi_1$, and the corresponding
$\mathbb{R}^*$-family of pseudospherical surfaces  is given by
$$F^\lambda=-\lambda\frac{\partial \Phi}{\partial
\lambda}\Phi_\lambda^{-1}:\mathbb{R}^2\to\mathbb{R}^3.$$

\section{Dressing Action}
Harmonicity equations for maps from a Riemann surface into a
compact symmetric space admit an infinite dimensional group of
symmetries:

\vspace{.20in}

Let $G$ be a compact semisimple Lie group and $\tau:G\rightarrow
G$ an involution with fixed set $K$. Fix an Iwasawa decomposition of $K^\C$: $K^\C=KB$, where $B$ is a solvable subgroup of $K^\C$. Fix $0<\varepsilon <1$. Let $C_\varepsilon$ and
$C_{1/\varepsilon}$ denote the circles of radius $\varepsilon$ and
$1/\varepsilon$ centered at $0\in\mathbb{C}$; define
$$I_\varepsilon=\left\{\lambda\in\mathbb{P}^1\,|\,\,|\lambda|<\varepsilon\right\},\,\,
I_{1/\varepsilon}=\left\{\lambda\in\mathbb{P}^1\,|\,\,|\lambda|>1/\varepsilon\right\},\,\,
E^\varepsilon=\left\{\lambda\in\mathbb{P}^1\,|\,\,\varepsilon<|\lambda|<1/\varepsilon\right\};$$
put $I^\varepsilon=I_\varepsilon\cup I_{1/\varepsilon}$ and
$C^\varepsilon=C_{\varepsilon}\cup C_{1/\varepsilon}$ so that
$\mathbb{P}^1=I^\varepsilon\cup C^\varepsilon\cup E^\varepsilon.$
Consider the infinite-dimensional twisted Lie groups
\begin{align*}
\Lambda^\varepsilon G^\C&=\big\{\gamma:C^\varepsilon\to
G^{\mathbb{C}}\,\,(\mathrm{smooth})\,|\,\,\mbox{$\tau\gamma(\lambda)=\gamma(-\lambda)$}\big\}\\
\Lambda_E^\varepsilon G^\C&=\left\{\gamma\in \Lambda^\varepsilon
G^\C \,|\,\,\mbox{$\gamma$ extends holomorphically to $\gamma:
E^\varepsilon\to G^\C$} \right\}\\
\Lambda_{I,B}^\varepsilon G^\C&=\left\{\gamma\in \Lambda^\varepsilon
G^\C \,|\,\,\mbox{$\gamma$ extends holomorphically to $\gamma:
I^\varepsilon\to G^\C$ and $\gamma(0)\in B$}\right\}.
\end{align*}

The basis of our action is the following decomposition:
\begin{thm}\cite{bG}\label{Fact}
\emph{The multiplication map
\begin{equation*}
\mu: \Lambda_E^\varepsilon G^\C \times \Lambda_{I,B}^\varepsilon
G^\C\rightarrow \Lambda^\varepsilon G^\C,\,\,\, (  g_E,g_I  )
\mapsto g_Eg_I
\end{equation*}
is a diffeomorphism onto an open dense subset ${\cal{U}}$ of a
component of $\Lambda^\varepsilon G^\C$.}
\end{thm}
\begin{rem}
Denote by $\Lambda^\varepsilon G$ the subgroup of
$\Lambda^\varepsilon G^\C$ formed by the loops $\gamma$ satisfying
the reality condition
$\gamma(\lambda)=\overline{\gamma(1/\overline{\lambda})}$. In this
case, the multiplication map $\mu$ gives a difeomorphism of
$\Lambda_E^\varepsilon G \times \Lambda_{I,B}^\varepsilon G$ onto
$\Lambda^\varepsilon G$. The dressing action on harmonic maps studied by Burstall and Pedit in \cite{BP2} is based on this decomposition of $\Lambda^\varepsilon G$.
\end{rem}

For $g_I\in \Lambda_{I,B}^\varepsilon G^\C$, let ${\cal{U}}^\varepsilon_{g_I}$ be
the open neighborhood of the identity $1$ in
$\Lambda_E^\varepsilon G^\C$ defined by: $g_E\in{\cal{U}}^\varepsilon_{g_I}$ if, and
 only if, there are unique $\hat{g}_E\in\Lambda_E^\varepsilon G^\C$ and $\hat{g}_I\in\Lambda_{I,B}^\varepsilon G^\C$ such that
$g_Ig_E=\hat{g}_E\hat{g}_I$
on $C^\varepsilon$. Write $g_I\#_\varepsilon  g_E$ for
$\hat{g}_E$. Thus $g_I\#_\varepsilon g_E=g_Ig_E\hat{g}_I^{-1}$.
One can prove easily the following:
\begin{lem}
\emph{a) ${\cal{U}}^\varepsilon_{1}=\Lambda_E^\varepsilon G$ and
$1\#_\varepsilon g_E= g_E$ for all $g_E\in \Lambda_E^\varepsilon
G$. b) For all $g_I\in \Lambda_{I,B}^\varepsilon G$,
$g_I\#_\varepsilon 1=1$. c) Let $g_1,g_2\in \Lambda_I^\varepsilon
G$, $g_E\in {\cal{U}}^\varepsilon_{g_1}$ and suppose $g_1\#_\varepsilon g_E\in
{\cal{U}}_{g_2}$ so that $g_2\#_\varepsilon(g_1\#_\varepsilon
g_E)$ is defined. Then, $g_E\in {\cal{U}}^\varepsilon_{g_2 g_1}$ and $(g_2
g_1)\#_\varepsilon g_E=g_2\#_\varepsilon(g_1\#_\varepsilon g_E)$.}
\end{lem}

Hence we conclude that $g_I\#_\varepsilon g_E$ defines a
local action of $\Lambda_I^\varepsilon G^\C$ on
$\Lambda_E^\varepsilon G^\C$.

 For $0<\varepsilon<\varepsilon'<1$ we have
injections $\Lambda_{I,B}^{\varepsilon'}G^\C\subset
\Lambda_{I,B}^{\varepsilon}G^\C$ and
$\Lambda_E^{\varepsilon}G^\C\subset \Lambda_E^{\varepsilon'}G^\C$.
Similarly, for $0<\varepsilon<1$, we have
$\Lambda_{\rm{hol}}G^\C\subset \Lambda_E^{\varepsilon}G^\C$, where
$$\Lambda_{\rm{hol}}G^\C=\bigcap_{0<\varepsilon <1}
\Lambda_E^\varepsilon G^\C.$$ Its is easy to see that
\begin{align*}
\Lambda_{\rm{hol}}G^\C & =\big\{\gamma:\C^*\to
G^\C\,|\,\,\mbox{$\gamma$ is holomorphic and
$\tau\gamma(\lambda)=\gamma(-\lambda)$} \big\}.
\end{align*}
The dressing actions are compatible with these inclusions:
\begin{thm}
\emph{For $0<\varepsilon<\varepsilon'<1$,
$\gamma\in\Lambda_{I,B}^{\varepsilon'}G^\C \subset
\Lambda_{I,B}^{\varepsilon}G^\C$, and $g\in\Lambda_E^\varepsilon
G^\C\subset\Lambda_E^{\varepsilon'}G^\C$, we have: $g\in {\mathcal{U}}_\gamma^{\varepsilon}$ if and only if $g\in {\mathcal{U}}_\gamma^{\varepsilon'}$;
$\gamma\#_{\varepsilon'} g=\gamma\#_\varepsilon g\in
\Lambda^\varepsilon_E G^\C.$}
\end{thm}
\begin{proof}
We argue as in \cite{BP2}, Proposition 2.3. Suppose that $g\in {\mathcal{U}}_\gamma^{\varepsilon'}$.
Then, on $C^{\varepsilon'}$ we can write
$$\gamma\#_{\varepsilon '}g=\gamma g(\gamma g)^{-1}_{I^{\varepsilon '}}.$$ Since
 $\gamma\#_{\varepsilon '}g$ has a holomorphic extension to $E^{\varepsilon '}$ while
 $\gamma g(\gamma g)^{-1}_{I^{\varepsilon '}}$ has an holomorphic extension to $I^{\varepsilon'}\cap E^\varepsilon$,
 it follows from a theorem of Painlev\'{e}
 that $\gamma\#_{\varepsilon '}g$ has a holomorphic extension to
 $E^{\varepsilon'}\cup C^{\varepsilon'}\cup(I^{\varepsilon'}\cap E^\varepsilon)=E^\varepsilon.$ Thus, $(\gamma g)_{I^{\varepsilon '}} \in\Lambda_{I,B}^\varepsilon G^\C$ while  $\gamma\#_{\varepsilon'}g\in \Lambda_E^\varepsilon G^\C$; which means that $g\in  {\mathcal{U}}_\gamma^{\varepsilon}$ and, from the uniqueness of the factorization of $\Lambda^\varepsilon G^\C$,  $\gamma\#_{\varepsilon}g=\gamma\#_{\varepsilon'} g$.

Conversely, suppose that $g\in {\mathcal{U}}_\gamma^{\varepsilon}$. Then, on $C^\varepsilon$ we can write
 $$(\gamma g)_{I^\varepsilon}=(\gamma g)^{-1}_{E^\varepsilon}\gamma g.$$ Since
$(\gamma g)_{I^\varepsilon}$ has a holomorphic extension to $I^{\varepsilon }$ while
 $(\gamma g)^{-1}_{E^\varepsilon}\gamma g$ has an holomorphic extension to $I^{\varepsilon'}\cap E^\varepsilon$,
  it follows from a theorem of Painlev\'{e} that $(\gamma g)_{I^\varepsilon}$  has a
  holomorphic extension to $I^{\varepsilon}\cup C^{\varepsilon}\cup(I^{\varepsilon'}\cap E^\varepsilon)=I^{\varepsilon'}.$
   Thus, $(\gamma g)_{I^{\varepsilon }} \in\Lambda_{I,B}^{\varepsilon'} G^\C$ while
    $\gamma\#_{\varepsilon}g\in \Lambda_E^{\varepsilon'} G^\C$; which means that
    $g\in  {\mathcal{U}}_\gamma^{\varepsilon'}$ and, from the uniqueness of the factorization of
    $\Lambda^{\varepsilon'} G^\C$,  $\gamma\#_{\varepsilon'}g=\gamma\#_{\varepsilon} g$.
\end{proof}

\begin{axiom}
\emph{The (local) action of each $\Lambda_{I,B}^{\varepsilon}G^\C$
preserves $\Lambda_{\rm{hol}}G^\C$ and, for
$0<\varepsilon<\varepsilon'<1,$
$\gamma\in\Lambda_{I,B}^{\varepsilon'}G^\C \subset
\Lambda_{I,B}^{\varepsilon}G$ and $g\in
\mathcal{U}_\gamma^\varepsilon\cap \Lambda_{\rm{hol}}G^\C$, we have $\gamma\#_{\varepsilon'}
g=\gamma\#_\varepsilon g$.}
\end{axiom}

 Henceforth, we write
$\gamma\# g$ for this (local) action on $\Lambda_{\rm{hol}}G^\C$.

Now, let $\Phi:\C\rightarrow  \Lambda_{\rm{hol}}G^\C$ be a smooth
map and $g_I\in \Lambda_{I,B}^\varepsilon G^\C$. Define the (smooth)
map $g_I\#
\Phi:\Phi^{-1}({\cal{U}}_{g_I}^\varepsilon)\subset\C\rightarrow
\Lambda_{\rm{hol}}G^\C$ by
\begin{equation*}
(g_I\# \Phi )(p)=g_I\#(\Phi(p))\,.
\end{equation*}
Observe that we can see a complex extended framing $\Phi$ as a map
into $\Lambda_{\rm{hol}}G^\C$, because $\Phi^{-1}{d}\Phi$ is holomorphic in
$\lambda$ on $\C^*$.

The relevance of the local action $\#$  is contained
in the following theorem:
\begin{thm}
\emph{If $\Phi: \C\rightarrow \Lambda_{\rm{hol}}G^\C$ is a complex extended
framing then so is $g_I\# \Phi$.}
\end{thm}
\begin{proof}
To see that $g_I\# \Phi$ is a complex extended
framing, write $g_I\Phi=ab$, where $a=g_I\# \Phi$ and
$b:\Phi^{-1}({\cal{U}}_{g_-})\subset \C\rightarrow \Lambda_{I,B}^\varepsilon
G^\C$. Then
\begin{equation}\label{lamhol}
a^{-1}{d}a=\mathrm{Ad}_b(\Phi^{-1}{d}\Phi-b^{-1}{d}b)\,,
\end{equation}
so that
\begin{equation}\label{lamhol1}
\lambda a^{-1}{d}a=\mathrm{Ad}_b(\lambda \Phi^{-1}{d}\Phi-\lambda
b^{-1}{d}b)
\end{equation}
and
\begin{equation}\label{lamhol2}
\lambda^{-1} a^{-1}{d}a={Ad}_b(\lambda^{-1}
\Phi^{-1}{d}\Phi-\lambda^{-1} b^{-1}{d}b)\,.
\end{equation}
Now, all the ingredients on the right side of (\ref{lamhol1}) are
holomorphic in $\lambda$ on a neighborhood of $0$ so that $\lambda
a^{-1}\mathrm{d}a$ is also; similarly, all the ingredients on the
right side of (\ref{lamhol2}) are holomorphic in  $\lambda$ on a
neighborhood of $\infty$ so that $\lambda^{-1} a^{-1}{d}a$ is
also. Whence, $a$ is a complex extended framing.
\end{proof}
From equation (\ref{lamhol}) we see that the $(1,0)$-part of
$(a^{-1}{d}a)_{\mathfrak{m}^\C}$ lies on an adjoint orbit of the
$(1,0)$-part of $(\Phi^{-1}{d}\Phi)_{\mathfrak{m}^\C}$. Then:
\begin{prop}\label{ncnc}
\emph{$\tilde{\varphi}^\lambda=a_\lambda\cdot \varphi_0$ is
conformal if, and only if, $\varphi^\lambda=\Phi_\lambda\cdot
\varphi_0$ is conformal.}
\end{prop}

\section{Simple factors}

Given $\gamma\in \Lambda_{I,B}^\varepsilon G^\C$ and $g\in
\Lambda_{\mathrm{hol}}G^\C$, a basic problem is to compute
$\gamma\# g$. This is a Riemann-Hilbert problem and, in general,
explicit solutions are not available. However, as the philosophy
underlying the work of Terng and Uhlenbeck \cite{Tu} suggests,
there should be certain elements of $\gamma\in \Lambda_{I,B}G^\C$,
the \textit{simple factors}, for which one can compute explicitly
$\gamma\# g$ by algebra alone. In this section we construct the
simple factors that are relevant to our geometric problem.

\vspace{.25in}

Let $L$ be an $1$-dimensional isotropic subspace of
$(\mathbb{R}^3)^\C\cong \C^3$: $(L,L)=0$. Let $Q\in\mathrm{SO}(3,\mathbb{R})$ be defined by (\ref{Q}),  suppose that $QL\neq L$
and consider the decomposition
\begin{equation*}
\C^3=L\oplus Q L\oplus L_0\,,
\end{equation*}
where $L_0=(L\oplus Q L)^\perp$. Denote by $\pi_L$, $\pi_{Q L}$
and $\pi_{L_0}$ the corresponding projections. For each
$\alpha\in\C\setminus\{0\}$ set
\begin{equation}\label{simple}
p_{\alpha,L}(\lambda)=\frac{\alpha-\lambda}{\alpha+\lambda}\pi_L+\pi_{L_0}+
\frac{\alpha+\lambda}{\alpha-\lambda}\pi_{Q L}\,.
\end{equation}
Thus, $p_{\alpha,L}:\mathbb{P}^1\setminus\{\pm\alpha\}\rightarrow
\mathrm{SO}(3,\mathbb{C})$ and  $p_{\alpha,L}(0)=\mathrm{1}$.
Moreover, each $p_{\alpha,L}$ is twisted, that is, $\tau
p_{\alpha,L}(\lambda)= p_{\alpha,L}(-\lambda)$. Then
$p_{\alpha,L}\in \Lambda_{I,B}^\varepsilon\mathrm{SO}(3,\mathbb{C})$
for some $\varepsilon <1$. The key to computing the
dressing action of $p_{\alpha,L}$ is the following proposition:
\begin{prop}\cite{Buiso}\label{simpledress}
\emph{Let $\Phi$ be a germ at $\alpha$ of a holomorphic map into
$\mathrm{SO}(3,\mathbb{C})$ such that $\tau
\Phi(\lambda)=\Phi(-\lambda)$. Suppose further that $Q
\Phi^{-1}(\alpha)L\neq \Phi^{-1}(\alpha)L$. Then
$p_{\alpha,\Phi^{-1}(\alpha)L}\in \Lambda_{I,B}^\varepsilon
{\mathrm{SO}}(3,\C)$  and
$$p_{\alpha,L}\Phi p^{-1}_{\alpha,\Phi^{-1}(\alpha)L}$$ is
holomorphic and invertible at $\alpha$.}
\end{prop}
\begin{axiom}
\emph{Let $<\!\!\varphi_0\!\!>$ be the subspace of $\mathbb{R}^3$
generated by $\varphi_0$ and denote by $<\!\!\varphi_0^\perp
\!\!>$ its real orthogonal complement in $\mathbb{R}^3$. Then,
given $g\in\Lambda_{\mathrm{hol}}\mathrm{SO}(3,\C)$,
$g\in{\cal{U}}^\varepsilon_{p_{\alpha ,L}}$ if, and only if,
$g^{-1}(\alpha )L$ is not contained in $<\!\!\varphi_0^\perp
\!\!>^\C$. For $g\in {\cal{U}}^\varepsilon_{p_{\alpha,L}}$, we
have
\begin{equation}\label{paL}
p_{\alpha,L}\# g=p_{\alpha,L}g p^{-1}_{\alpha,g^{-1}(\alpha)L}\,.
\end{equation}}
\end{axiom}
\begin{proof}
The eigenspaces of $Q$ are $<\!\!\varphi_0\!\!>$ and
$<\!\!\varphi_0^\perp \!\!>$. Hence, since $g^{-1}(\alpha)L$ is
isotropic and no real subspace is isotropic, $Q
g^{-1}(\alpha)L\neq g^{-1}(\alpha)L$ if, and only if,
$g^{-1}(\alpha)L$ is not contained in $<\!\!\varphi_0^\perp
\!\!>^\C$.

The first part of Proposition \ref{simpledress} assures us that
$p_{\alpha,g^{-1}(\alpha)L}\in \Lambda_{I,B}^\varepsilon
\mathrm{SO}(3,\mathbb{C})$. So we just need to prove that
$p_{\alpha,L}\# g$ given by (\ref{paL}) is an element of
$\Lambda_{\mathrm{hol}}\mathrm{SO}(3,\mathbb{C})$. Clearly
$p_{\alpha,L}\# g$ is twisted, since it is a product of maps with
this property. The holomorphicity at $\alpha$ follows directly
from Proposition \ref{simpledress} and then we get holomorphicity
at $\pm \alpha$ from the twisting condition.
\end{proof}

\begin{rem}
In \cite{Dol}, the authors gave the following general definition
of simple factors: Let $G$ be a compact Lie group and
$\rho:G^\C\to {\mathrm{GL}}(V)$
 a representation. A semisimple element $H\in \mathfrak{g}^\C$ is
said to be \emph{$\rho$-integral} if $\rho(H)\in{\mathrm{End}}(V)$
has only integer values. For any $\alpha\in \C\setminus\mathbb{R}$
and $\rho$-integral element $H\in i\mathfrak{g}$, the loop
\begin{equation}\label{D}
p_{\alpha,H}(\lambda)=\exp\Big(\ln\Big(\frac{\lambda-\alpha}{\lambda-\overline{\alpha}}\Big)
H\Big)
\end{equation}
 is a \emph{simple factor}, with $\lambda$ belonging to $\mathbb{R}^*$, the multiplicative group of nonzero real numbers.
The condition $H\in i\mathfrak{g}$ ensures that $p_{\alpha,H}$
satisfies the reality condition
$p_{\alpha,H}(\lambda)=\overline{p_{\alpha,H}(\overline{\lambda})}$.
Hence, on the real axis $p_{\alpha,H}$ takes values in $G$.
Moreover, if $N=G/K$ is a symmetric space with automorphism $\tau$
and associated symmetric decomposition $\g=\lk\oplus\mathfrak{m}$,
then the conditions $H\in i\mathfrak{m}$ and $\alpha=ir$, with
$r\in\mathbb{R}\setminus \{0\}$, ensure that $p_{\alpha,H}$ is
also twisted: $\tau p_{\alpha,H}(\lambda)=
p_{\alpha,H}(-\lambda)$. In our case, we are not imposing any
reality condition to the simple factors (\ref{simple}), which can
also be written as:
$$p_{\alpha,H}(\lambda)=\exp\Big(\ln\Big(\frac{\alpha-\lambda}{\alpha+\lambda}\Big) H\Big),$$
where $H\in\mathfrak{m}^\C$ is $\rho$-integral with respect to the
standard representation of $\mathrm{SO}(3,\C)$.
\end{rem}

\begin{rem}\label{table} Let $L$ be an $1$-dimensional isotropic subspace of
$(\mathbb{R}^3)^\C\cong \C^3$. Suppose that $QL\neq L$. The cross
product multiplication table for the decomposition $\C^3=L\oplus Q
L\oplus L_0$, where $L_0=(L\oplus Q L)^\perp$, is given by

\begin{center}
\small{\begin{tabular}{|c||c|c|c|}
  \hline
  $\times$ & $L$ & $L_0$ & $QL$ \\
  \hline\hline
  $L$ & $0$ & $L$ & $L_0$ \\\hline
  $L_0$ & $L$ & $0$ & $QL$ \\\hline
  $QL$ & $L_0$ & $QL$ & 0 \\
  \hline
\end{tabular}}
\end{center}
\end{rem}

\section{Bianchi-B\"{a}cklund transforms via dressing actions}

In this section we prove that the actions of these simple factors
amount to Bianchi-B\"{a}cklund transformations.

\vspace{.25in}

Start with an everywhere non-conformal harmonic map
$\varphi:\C\rightarrow S^2$. Let $\Phi:\C\rightarrow
\Lambda_{\mathrm{hol}}{\mathrm{SO}}(3,\C)$ be an extended framing
associated to $\varphi$. By applying  formula (\ref{adlersymes})
to $\Phi$, we get a map $F^\lambda:\C\rightarrow \mathbb{R}^3$,
for each $\lambda\in \C^*$, such that
${d}F^\lambda=[\varphi^\lambda,*{d}\varphi^\lambda ]$, that is, a
CGC $K=1$ surface (which is real when $\lambda\in S^1$) without
umbilics, with normal $\varphi^\lambda=\Phi_\lambda\cdot
\varphi_0$. Assume that $\Phi_{\lambda}(z_0)=\mathrm{1}$ for all
$\lambda\in \C^*$.

Choose $\alpha\in\C\setminus\{0,\pm \lambda\}$. Consider  the
action of a simple factor $p_{\alpha,L}$ on $\Phi$:
\begin{equation*}
\tilde{\Phi}=p_{\alpha,L}\# \Phi= p_{\alpha,L}\Phi
p_{\alpha,\tilde{L}}^{-1}:\Phi^{-1}\big({\mathcal{U}}^\varepsilon_{p_{\alpha,L}}\big)\to
\Lambda_{\mathrm{hol}}\mathrm{SO}(3,\C)\,,
\end{equation*}
where $\tilde{L}=\Phi(\alpha)^{-1}L$. Set $h=p_{\alpha,L}$ and
$\tilde{h}=p_{\alpha,\tilde{L}}$. Since
\begin{align*}
-i\lambda \frac{\partial \tilde{\Phi}}{\partial
\lambda}&\tilde{\Phi}_\lambda^{-1}= -i\lambda
h_\lambda\Big(\frac{\partial \Phi}{\partial
\lambda}\Phi_\lambda^{-1}
-\Phi_\lambda\tilde{h}_\lambda^{-1}\frac{\partial
 \tilde{h}}{\partial
\lambda}\Phi_\lambda^{-1}+h_\lambda^{-1}\frac{\partial h}{\partial
\lambda}\Big)h_\lambda^{-1}\,,
\end{align*}
our new  CGC $K=1$ surface  equals
\begin{equation}\label{drecgc}
\tilde{F}^\lambda=F^\lambda+i\lambda\Phi_\lambda\tilde{h}_\lambda^{-1}\frac{\partial
 \tilde{h}}{\partial
\lambda}\Phi_\lambda^{-1}\,,
\end{equation}
up to  (complex) Euclidean motions -- translation by
$h_\lambda^{-1}\frac{\partial h}{\partial \lambda}$ followed by a
rotation (conjugation by $h_\lambda$). The corresponding normal is
$\tilde{\varphi}^\lambda=\Phi_\lambda \tilde{h}_\lambda^{-1}\cdot
\varphi_0.$ By Proposition \ref{ncnc}, $\tilde{\varphi}^\lambda$
is also everywhere non-conformal, whence $\tilde{F}^\lambda$ has
no umbilic points.

\begin{thm}\label{bbdress}
\emph{$\tilde{F}^\lambda\in{\mathbf{BB}}_{\beta}(F^\lambda)$, with
$\beta=\ln \big(\alpha/\lambda\big)$. Moreover, any
Bianchi-B\"{a}cklund transform of $F$ is of the form
(\ref{drecgc}) for some simple factor $p_{\alpha,L}$.}
\end{thm}
\begin{proof}
1) Since $\varphi^\lambda$ and $\tilde{\varphi}^\lambda$ are harmonic with respect
to the same conformal structure $z=x+iy$ on $\C$, we have
$\mathrm{\Pi}_{F^\lambda}\big(\frac{\partial}{\partial
z},\frac{\partial}{\partial
z}\big)=\mathrm{\Pi}_{\tilde{F}^\lambda}\big(\frac{\partial}{\partial
z},\frac{\partial}{\partial z}\big)=0$, that is, $z$ is a
curvature line coordinate with respect to $F^\lambda$ and $\tilde{F}^\lambda$.

2)  Since
\begin{equation*}
\tilde{h}_\lambda^{-1}\frac{\partial
 \tilde{h}}{\partial
\lambda}=A_\lambda(\alpha)\big(\pi_{Q\tilde{L}}-\pi_{\tilde{L}}\big)\,,
\end{equation*}
where $A_\lambda(\alpha)=\frac{2\alpha}{(\alpha-\lambda)(\alpha+\lambda)},$ we have:
\begin{equation}\label{fs}
(\tilde{F}^\lambda-F^\lambda,\tilde{F}^\lambda-F^\lambda) =
\frac{\lambda^2}{2}\mathrm{Tr}\,\big\{A_\lambda(\alpha)\Phi_\lambda\big(\pi_{Q\tilde{L}}-
\pi_{\tilde{L}})\Phi_\lambda^{-1}\big\}^2 =
\lambda^2A_\lambda^2(\alpha)\,,
\end{equation}
that is, $\tilde{F}^\lambda-F^\lambda$ has constant length.

3) Since $[Q,\varphi_0]=0$ and
$\pi_{Q\tilde{L}}=Q\pi_{\tilde{L}}Q^{-1}$, one can easily check
that
\begin{eqnarray*}
(\tilde{F}^\lambda-F^\lambda,\varphi^\lambda)=(\tilde{F}^\lambda-F^\lambda,\tilde{\varphi}^\lambda)=0.
\end{eqnarray*}

4) Recall that $\varphi_0=\xi_{e_1}$, where $e_1,e_2,e_3$ is the
canonical basis of $\mathbb{R}^3$ and $\xi_{e_1}(v)=e_1\times v$.
Since $\tilde{L}$ is isotropic and $Q\tilde{L}\neq\tilde{L}$,
$\tilde{L}$ is generated by a vector $\tilde{v}$ of the form
$\tilde{v}=\frac{1}{2}e_1+\tilde{a}e_2+\tilde{b}e_3$, with
$\tilde{a}^2+\tilde{b}^2=-\frac{1}{4}$. Note that
$e_1=\tilde{v}+Q\tilde{v}$.

Now:
\begin{equation*}
(\varphi^\lambda,\tilde{\varphi}^\lambda)=-\frac{1}{2}\mathrm{Tr}\,
\big\{\Phi_\lambda\varphi_0\Phi_\lambda^{-1}\Phi_\lambda\tilde{h}_\lambda^{-1}\varphi_0\tilde{h}_\lambda
\Phi_\lambda^{-1}\big\}=-\frac{1}{2}\mathrm{Tr}\,
\big\{\underbrace{\varphi_0\tilde{h}_\lambda^{-1}\varphi_0\tilde{h}_\lambda}_{\equiv\rho}\big\}\,.
\end{equation*}
  Set $a_\lambda(\alpha)=\frac{\alpha+\lambda}{\alpha-\lambda}$.
  Fix $X\in\tilde{L}_0$, $Y\in\tilde{L}$ and $Z\in Q\tilde{L}$.
 By using the cross product multiplication table with respect to
the decomposition  $\C^3=\tilde{L}\oplus \tilde{L}_0\oplus
Q\tilde{L}$ (see Remark
 \ref{table}) and the  triple product expansion $\vec{a}\times (\vec{b}\times
 \vec{c})=(\vec{a},\vec{c})\vec{b}-\vec{c}(\vec{a},\vec{b})$, one
 can straightforwardly compute $\rho(X)$, $\rho(Y)$ and $\rho(Z)$.
 For example:
\begin{eqnarray*}
\rho(X)&\!\!\!\!=\!\!\!\!&\varphi_0\tilde{h}_\lambda^{-1}\varphi_0\tilde{h}_\lambda(X)=\varphi_0\tilde{h}_\lambda^{-1}
\big(\underbrace{\tilde{v}\times X}_{\in \tilde{L}}+
\underbrace{Q\tilde{v}\times X}_{\in Q\tilde{L}}\big)
=\varphi_0\big(a_\lambda(\alpha) \tilde{v}\times
X+a_\lambda(\alpha)^{-1}Q\tilde{v}\times
X\big)\\&\!\!\!\!=\!\!\!\!&
a_\lambda(\alpha)Q\tilde{v}\times\big(\tilde{v}\times
X\big)+a_\lambda(\alpha)^{-1}\tilde{v}\times\big(Q\tilde{v}\times
X\big)= -\frac{a_\lambda(\alpha)+a_\lambda(\alpha)^{-1}}{2}X.
\end{eqnarray*}
Similarly, we have:
$$\rho(Y)=-\frac{a_\lambda(\alpha)^{-1}}{2}Y+a_\lambda(\alpha)^{-1}\underbrace{Q\tilde{v}\times
\big(Q\tilde{v}\times Y\big)}_{\in Q\tilde{L}}$$ and
$$\rho(Z)=-\frac{a_\lambda(\alpha)}{2}Z+a_\lambda(\alpha)\underbrace{\tilde{v}\times
(\tilde{v}\times Z)}_{\in \tilde{L}}.$$ Then,
\begin{equation*}
\cos\sigma=(\varphi^\lambda,\tilde{\varphi}^\lambda)=-\frac{1}{2}\mathrm{Tr}\,\rho=
\frac{a_\lambda(\alpha)+a_\lambda(\alpha)^{-1}}{2},\,
\end{equation*}
and we conclude that $(\varphi^\lambda,\tilde{\varphi}^\lambda)$ is constant.

\vspace{.10in}

So, $\tilde{F}^\lambda$ is a Bianchi-B\"{a}cklund transformation of $F^\lambda$.
Next we shall find the corresponding $\beta$ parameter.

By (\ref{newcgc}) and (\ref{fs}), we get $$\frac{1}{\sinh
\beta}=\pm \lambda A_\lambda(\alpha).$$ Taking account  (\ref{osigma}) and
(\ref{newcgc}), one can compute
 $\sin\sigma$:
 \begin{equation*}
\sin\sigma=\big(\varphi^\lambda\times\tilde{\varphi}^\lambda,\sinh{\beta}(F^\lambda-\tilde{F}^\lambda)\big)=\pm
i\big(\varphi^\lambda\times\tilde{\varphi}^\lambda,\Phi_\lambda(\pi_{Q\tilde{L}}-\pi_{\tilde{L}})\Phi_\lambda^{-1}\big)=
\pm i\frac{a_\lambda(\alpha)-a_\lambda(\alpha)^{-1}}{2}\,.
\end{equation*}
Hence,
\begin{equation*}
-i\cosh\beta=\cot\sigma=\mp
i\frac{a(\alpha)+a(\alpha)^{-1}}{a(\alpha)-a(\alpha)^{-1}}\,.
\end{equation*}
By Lemma \ref{-beta}, we can take
$$\frac{1}{\sinh\beta}=\lambda A_\lambda(\alpha)=\frac{a_\lambda(\alpha)-a_\lambda(\alpha)^{-1}}{2}\,\,\,\,\,\,\mbox{and}\,\,\,\,\,\,
\cosh{\beta}=\frac{a_\lambda(\alpha)+a_\lambda(\alpha)^{-1}}{a_\lambda(\alpha)-a_\lambda(\alpha)^{-1}};$$
hence $\beta=\ln(\alpha/\lambda)$.

\vspace{.10in}

It remains to prove the converse, that is, any
Bianchi-B\"{a}cklund transformation of $F^\lambda$ amounts to the
  dressing action of some simple factor $p_{\alpha,L}$  (up to conjugation by a constant complex
  matrix):

  Evaluating $\tilde{F}^\lambda-F^\lambda$ at $z_0$ gives
\begin{equation*}
\tilde{F}^\lambda(z_0)-F^\lambda(z_0)=i\lambda
A_\lambda(\alpha)(\pi_{QL}-\pi_L)\,.
\end{equation*}
For some fixed spectral parameter $\beta\in\C\setminus\{in\pi,
n\in\mathbb{Z}\}$, we choose $\alpha=e^{\beta}\lambda$.
 On the other hand, if $L$ is
 generated by $v=\frac{1}{2}e_1+ae_2+be_3$, for some
$a,b\in\C$ such that $a^2+b^2=-\frac{1}{4}$, we have
\begin{equation*}
\pi_{Q{L}}-\pi_{{L}}=
  \begin{pmatrix}
    0 & 2{a} & 2{b} \\
     -2 {a} & 0 & 0 \\
-2{b} & 0 & 0
  \end{pmatrix}\,,
  \end{equation*}
  and so any non-isotropic direction in $<\!\!\varphi_0^\perp\!\!>^\C$
  is generated by a vector of the form
  $\pi_{Q{L}}-\pi_{{L}}$. We therefore deduce from the
  uniqueness of solutions to the Bianchi-B\"{a}cklund PDEs that
  each Bianchi-B\"{a}cklund transformation of $F^\lambda$ amounts to the
  dressing action of some simple factor $p_{\alpha,L}$ (up to conjugation by a constant complex
  matrix).
\end{proof}

In the pseudospherical case, Lie observed that every B\"{a}cklund
tranformation is a combination of transformations of Lie
 and Bianchi. In the spherical case,
the tangent planes at corresponding points on $F$ and
$\tilde{F}=\mathbf{BB}_{i\pi/2}(F)$ are orthogonal, and we have:
\begin{axiom}
\emph{$\mathbf{BB}_\beta= \mathbf{S}_{\beta+i\pi/2}^{-1}\circ
\mathbf{BB}_{i\pi/2}\circ \mathbf{S}_{\beta+i\pi/2}.$}
\end{axiom}

\vspace{.20in}

\textit{B\"{a}cklund transforms and  dressing actions.} In
\cite{U} (see also \cite{GS}), Uhlenbeck made the observation
 that B\"{a}cklund transforms for pseudospherical surfaces (CGC $K<0$  surfaces)
 amount to dressings of the simplest type. She used $\mathrm{SU(2)}$ as  symmetry group.
 If we want to use  $\mathrm{SO}(3,\R)$
as symmetry group, we shall proceed as follows:

The basis of the action is the same decomposition of Theorem
\ref{Fact}, restricted to the loops satisfying the reality
condition $\overline{\xi(\lambda)}=\xi(\bar{\lambda})$, and the
simple factors are those of the form (\ref{D}) with $H\in
i\mathfrak{m}$ and $\alpha\in i\mathbb{R}$. More explicitly, the
simple factors we need in this case are of the form
\begin{equation}\label{boing}
p_{\alpha,L}(\lambda)=\frac{\alpha-\lambda}{\alpha+\lambda}\pi_L+\pi_{L_0}+
\frac{\alpha+\lambda}{\alpha-\lambda}\pi_{{QL}},
\end{equation}
with $QL=\overline{L}$.

\begin{rem}
The simple factors (\ref{boing}) generate the group of all
rational maps $\gamma: \mathbb{P}^1\to \mathrm{SO}(3,\C)$,
holomorphic at $0$ and $\infty$,      satisfying the reality
condition $\overline{\gamma(\lambda)}=\gamma(\bar{\lambda})$ and
the twisted condition $\tau\gamma(\lambda)=\gamma(-\lambda)$ (cf.
\cite{Dol}).
\end{rem}

\subsection{Bianchi-B\"{a}cklund Permutability theorem }

Taking account the results we have obtained above, the
Bianchi-B\"{a}cklund Permutability theorem is a direct consequence
of the following:
\begin{thm}
\emph{Consider simple factors $p_{\alpha_1,L_1}$ and
$p_{\alpha_2,L_2}$ with $\alpha_1^2\neq\alpha_2^2$. Set
$L'_1=p_{\alpha_2,L_2}(\alpha_1)L_1$ and
$L'_2=p_{\alpha_1,L_1}(\alpha_2)L_2$. Assume that $QL'_i\neq
L'_i$, $i=1,2$. Then
\begin{equation}\label{bbp}
p_{\alpha_1,L'_1}p_{\alpha_2,L_2}=p_{\alpha_2,L'_2}p_{\alpha_1,L_1}.\end{equation}}
\end{thm}
This theorem is a simple adaptation of Proposition 4.15 of
\cite{Buiso} to our setting and its proof can be carried out
similarly, without any reality assumption.

\subsection{Getting a real solution from an old one}

If we want to obtain a new real CGC $K=1$ surface from an old one,
we have to perform two dressing actions with simple factors, as
the classical theory suggests.

\vspace{.25in}

For each pair $(\alpha,L)$, we introduce the holomorphic map
$q_{\alpha,L}:\mathbb{P}^1\setminus\{0,\pm\alpha\}\rightarrow
\mathrm{SO}(3,\C)$ defined by
\begin{equation*}
q_{\alpha,L}(\lambda)=p_{\alpha ,L}(\infty)p_{\alpha
,L}(\lambda)=\frac{\lambda-\alpha}{\lambda+\alpha}\pi_L+\pi_{L_0}+
\frac{\lambda+\alpha}{\lambda-\alpha}\pi_{QL}\,,
\end{equation*}
 where, as before,  $L$ is an
isotropic line in $(\mathbb{R}^3)^\C\cong \C^3$ such that $QL\neq
L$. Consider the automorphism $\mathbf{{R}}:\Lambda
\mathrm{SO(3,\C)}\to \Lambda\mathrm{SO(3,\C)}$ defined by
$\mathbf{{R}}(\gamma)(\lambda)=\overline{\gamma(1/\overline{\lambda})}.$
Clearly, $\mathbf{{R}}$ is an involution.  One can easily check
that
\begin{equation}\label{rsf}
{\mathbf{{R}}}(p_{\alpha,L})(\lambda)=
q_{\frac{1}{\overline{\alpha}},\overline{L}}(\lambda)\,.
\end{equation}
\begin{thm}\label{realth}
\emph{Set $L_1=L$, $L_2=\overline{L}$, $\alpha_1=\alpha$,
$\alpha_2=\frac{1}{\overline{\alpha_1}}$,
$L'_2=p_{\alpha_1,L_1}(\alpha_2)L_2$ and
$L'_1=p_{\alpha_2,L_2}(\alpha_1)L_1$. There exists $k\in
{K^\tau}^\C$ such that $kp_{\alpha_2,L'_2}p_{\alpha_1,L_1}\in
\Lambda_{I,B}^\varepsilon \mathrm{SO}(3,\mathbb{R})$ for some
$0<\varepsilon<1$.}
\end{thm}
\begin{proof}
Set $u=p_{\alpha_2,L'_2}p_{\alpha_1,L_1}$. We want to find
$k\in{K^\tau}^\C$ satisfying $\mathbf{R}(ku)=ku;$ and this is the
same to $k^{-1}\bar{k}=u\mathbf{R}(u)^{-1}.$ Now, it follows from
(\ref{bbp}) and (\ref{rsf}) that:
\begin{eqnarray}\label{pa}
u\mathbf{R}(u)^{-1}
&=&p_{\alpha_2,L'_2}p_{\alpha_1,L_1}q_{\alpha_2,L_2}^{-1}
q_{\alpha_1,\overline{L'_2}}^{-1}=p_{\alpha_1,L'_1}p_{\alpha_2,L_2}q_{\alpha_2,L_2}^{-1}
q_{\alpha_1,\overline{L'_2}}^{-1} \nonumber \\&=&
p_{\alpha_1,L'_1}p_{\alpha_2,L_2}(\infty)q_{\alpha_1,\overline{L'_2}}^{-1}.
\end{eqnarray}
But
\begin{equation*}\label{qqla}
\overline{L'_2}=\overline{p_{\alpha_1,L_1}
(\alpha_2)L_2}=q_{\alpha_2,L_2}(\alpha_1)L_1
=p_{\alpha_2,L_2}(\infty)L'_1 \,.
\end{equation*}
Hence
\begin{equation}\label{ns}
q_{\alpha_1,\overline{L'_2}}=p_{\alpha_2,L_2}(\infty)q_{\alpha_1,{L'_1}}p_{\alpha_2,L_2}(\infty).
\end{equation}
From (\ref{pa}) and (\ref{ns}), we obtain
$$u\mathbf{R}(u)^{-1}=p_{\alpha_1,L'_1}(\infty)p_{\alpha_2,L_2}(\infty).$$
Then, $P=u\mathbf{R}(u)^{-1}$ does not depend on $\lambda$ and,
since all the simple factors are twisted, we have
$P\in{K^\tau}^\C$. Finally, observe that
$$P\overline{P}=P\mathbf{R}(P)=u\mathbf{R}(u)^{-1}\mathbf{R}(u\mathbf{R}(u)^{-1})=1.$$
The existence of such $k\in{K^\tau}^\C$ follows now from Lemma
\ref{util}.
\end{proof}

\begin{lem}\label{util}
\emph{If $P\in {K^\tau}^\C$ satisfies $P\overline{P}=1$, then
$P=k^{-1}\overline{k}$ for some $k\in {K^\tau}^\C$.}
\end{lem}
\begin{proof}
The complex group $\C^*$ double covers ${K^\tau}^\C$ via
\begin{equation*}
z\in \C^*\mapsto \rho(z)=\mathrm{Ad}{\begin{pmatrix}
  z & 0 \\
  0 & z^{-1}
\end{pmatrix}}\in
\mathrm{Ad}_{\mathrm{SL}(2,\C)}\cong\mathrm{SO(3,\C)}\,.
\end{equation*}
Observe that $\rho(\overline{z})=\overline{\rho(z)}^{-1}$. Given
$P\in {K^\tau}^\C$ such that $P\overline{P}=1$, one can find
$z_0\in \C^*$ such that $\rho(z_0)=P$ and
$z_0\overline{z_0}^{-1}=1$ ; in particular, we can fix $z_0\in
\mathbb{R}^+$. This means that there exists $u_0\in \C^*$ such
that $z_0=u_0\overline{u_0}$. Hence
\begin{equation*}
P=\rho(z_0)=\rho(u_0)\rho(\overline{u_0})=\rho(u_0)\overline{\rho(u_0)}^{-1}.
\end{equation*}
Set $k=\rho(u_0)^{-1}$. Then $P=k^{-1}\overline{k}$.
\end{proof}

Let $\Phi:\C\rightarrow \Lambda_{\mathrm{hol}}
\mathrm{SO}(3,\mathbb R)$ be an extended framing associated to an
everywhere non-conformal harmonic map $\varphi:\C\rightarrow S^2$.
Applying  formula (\ref{adlersymes}) to $\Phi$, and evaluating at
$\lambda=1$, we get a real CGC $K=1$ surface $F$. With the
notations of Theorem \ref{realth},
\begin{equation*}
\tilde{\Phi}=(kp_{\alpha_2,L'_2}p_{\alpha_1,L_1})\# \Phi
\end{equation*}
is a new (real) extended framing.  Therefore, formula
(\ref{adlersymes}) applied to $\tilde{\Phi}$ gives a new real CGC
$K=1$ surface $F^*$, evaluating again at $\lambda=1$. Up to an
Euclidean motion, this surface is obtained out of $F$ by applying
two successive Bianchi-B\"{a}cklund transformations: if
$\beta_1\in\C\setminus\{in\pi, n\in \mathbb{Z}\}$ is such that
$\beta_1=\ln\alpha_1$, then, taking account Lemma \ref{-beta} and
Theorem \ref{bbdress},  $F^*$ belongs to
${\mathbf{BB}}_{i\pi-\overline{\beta}_1}\big({\mathbf{BB}}_{\beta_1}(F)\big),$
up to an Euclidean motion, which agrees with Theorem
\ref{realcondb}.

\begin{rem}
 The multiplication map $\mu$, restricted to loops $\gamma$ satisfying
the reality condition
$\gamma(\lambda)=\overline{\gamma(1/\overline{\lambda})}$, gives a
difeomorphism of $\Lambda_E^\varepsilon \mathrm{SO}(3,\mathbb{R})
\times \Lambda_{I,B}^\varepsilon \mathrm{SO}(3,\mathbb{R})$ onto
$\Lambda^\varepsilon \mathrm{SO}(3,\mathbb{R})$. Hence, the
extended framing
$\tilde{\Phi}=(kp_{\alpha_2,L'_2}p_{\alpha_1,L_1})\# \Phi$ is
well-defined everywhere.
\end{rem}

\end{document}